\documentclass[11pt]{amsart}
\usepackage[latin1]{inputenc}
\usepackage{eucal,bbm,mathrsfs}
\theoremstyle{plain}    
\newtheorem{thm}{Theorem}
\newtheorem*{thm*}{Theorem}
\theoremstyle{plain}    
\newtheorem{lem}{Lemma} 
\theoremstyle{plain}    
\newtheorem{prop}{Proposition} 
\theoremstyle{definition}
\newtheorem{defn}{Definition}
\theoremstyle{definition}
 
\theoremstyle{remark}
\newtheorem{rem}{Remark}

\newcommand{\ric}{\mathrm{Ric}}
\newcommand{\bundle}[1]{\CMcal{#1}}
\newcommand{\pen}{\#}
\newcommand{\curv}[1]{\mathcal{#1}}
\newcommand{\hodge}{\boldsymbol \star}
\newcommand{\schout}{\frak{k}}
\newcommand{\snabla}{\text {\tiny{$\! \nabla  $}}}

\begin{document}

\title{Conformally invariant Cotton and Bach tensor in $N$--dimensions}
\author{Mario Listing}
\address{Department of Mathematics, Stony Brook University, Stony Brook, NY 11794-3651, USA}
\email{listing@math.sunysb.edu}
\thanks{Supported by the German Research Foundation}
\begin{abstract}
This paper presents conformal invariants for Riemannian manifolds of dimension greater than or equal to four whose vanishing is necessary for a Riemannian manifold to be conformally related to an Einstein space. One of the invariants is a modification of the Cotton tensor, the other is a $n$--dimensional version of the Bach tensor. In general both tensors are smooth only on an open and dense subset of $M$, but this subset is invariant under conformal transformations. Moreover, we generalize the main result of "Conformal Einstein Spaces in $N$--Dimensions" published in \emph{Ann.~Global Anal.~Geom.} {\bf 20}(2) (2001).
\end{abstract}
\keywords{conformal invariants, Weyl tensor, Einstein space} 
\subjclass[2000]{Primary 53A30, Secondary 53A55}
\maketitle
\section{Introduction}
Let $(M^n,g)$ be a Riemannian manifold and $A$ be a tensor determined by $g$. $A$ is said to be \emph{conformally invariant} if for any local conformal transformation $\overline{g}=\psi ^{-2}g$, the corresponding tensors are related by $\overline{A}=\psi ^{2k}A$ with $k\in \mathbbm{Z}$. $k$ is called the \emph{weight} of $A$. In \cite{FeGr1} Fefferman and Graham considered scalar conformal invariants. They claimed that in the odd--dimensional case, each scalar conformal invariant is a linear combination of conformal Weyl invariants. Non--trivial conformal invariants of higher degree are the \emph{Cotton} tensor in dimension $n=3$, the \emph{Bach} tensor in dimension $n=4$ and the \emph{Weyl} tensor in dimension $n\geq 4$. 

The aim of the present paper is to find conformal invariants being necessary respectively sufficient for Riemannian manifolds to be locally conformally related to Einstein spaces. Because of the interest in physics, this problem is often studied in four dimensional Lorentz geometry (cf.~\cite{KNT,KNN,NP,Wunsch,CzMcLW}). The difficulty in the non--Riemannian case are the conformal transformations when the conformal factor has light--like gradient. In \cite{List3}, we found necessary and sufficient conditions for a $n$--dimensional semi--Riemannian manifold with $\det \curv{W}\neq 0$ to be conformally related to an Einstein space, where $\curv{W}$ is the Weyl tensor considered as endomorphism on $\Lambda ^2(T^*M)$. Moreover, we solved the problem (even if $\det \curv{W}=0$) in the four dimensional Riemannian case. In this paper, we extend the last result to $n$--dimensions and prove that the involved tensors are conformal invariants.

If $(M,g)$ is a Riemannian manifold of dimension $n\geq 4$ and $W$ the Weyl tensor, the set
\[
\bundle{E}:=\{ v\in TM\ |\  W(v,.,.,.)=0\}
\]
is conformally invariant and a vector space in each fiber. Let $M_\bundle{E}$ be the set of all points in $M$ at which the rank of $\bundle{E}$ is locally constant, then $M_\bundle{E}$ is open and dense in $M$. In section 3 we define the (smooth) vector field $\mathbb{T}:M_\bundle{E}\to TM_\bundle{E}$ by
\[
g(\mathbb{T},.):=\frac{1}{n-3}\frak{w}^\pen \left( \sum _ie_i\llcorner \curv{W}(\delta W_{e_i})\right) ,
\] 
where $\frak{w}^\pen :T^*M_\bundle{E}\to T^*M_\bundle{E}$ is the Moore--Penrose inverse of the negative Ricci contraction of $\curv{W}^2$ and $e_1,\ldots ,e_n$ is an orthonormal base of $T_pM$. The $(0,3)$ tensor
\[
C_{\mathbb{T}}:=\frac{1}{n-2}\ \mathrm{d} ^\snabla \left( \mathrm{Ric}-\frac{\mathrm{S}}{2(n-1)} g\right)-W(.,.,.,\mathbb{T})
\]
is a well defined conformally invariant of zero weight on $(M_\bundle{E},g)$. Obviously, $C_{\mathbb{T}}$ is the Cotton tensor if $\mathbb{T}=0$ or $n=3$. If $(M,g)$ is locally conformally related to a space with harmonic Weyl tensor (e.g.~if $g$ is conformally Einstein), $C_{\mathbb{T}}$ vanishes.

\begin{thm*}
Suppose $(M,g)$ is a simply connected Riemannian manifold of dimension $n\geq 4$ and with $\mathrm{rank}(\bundle{E})=0$ on $M_\bundle{E}$. Then $g$ is (globally) conformally related to an Einstein space if and only if $\mathbb{T}$ is extendible to a vector field on $M$ and
\[
E_{\mathbb{T}}:=\mathrm{Ric}-\frac{S}{n}g+(n-2)\left( \nabla \mathbb{T}^*+\mathbb{T}^*\otimes \mathbb{T}^* -\frac{1}{n}\left(\mathrm{div}(\mathbb{T})+|\mathbb{T}|^2\right)g\right)
\]
vanishes.
\end{thm*}
This is a major extension of the previous result in \cite{List3}, since a non--degenerate two form in the image of $\curv{W}_p$ implies $\mathrm{rank}(\bundle{E}_p)=0$, in particular if $n$ is even, $\mathrm{Im}(\curv{W})$ could be two dimensional and $\mathrm{rank}(\bundle{E})=0$. If $\mathrm{rank}(\bundle{E})$ is non--zero in some component of $M_\bundle{E}$, $E_{\mathbb{T}}=0$ is still sufficient for a conformal Einstein space, but it is no longer necessary. 

In the last section we show that the symmetric $(0,2)$ tensor
\[
\begin{split}
B_{\mathbb{T}}:=\delta _1\delta _4W+\frac{n-3}{n-2}\curv{W}(\mathrm{Ric})-\frac{n-4}{n-2}\curv{W}(\mathrm{sym}(E_{\mathbb{T}}))\qquad \qquad \quad \\
-(n-3)(n-4)\bigl[ \curv{W}(\mathbb{T}^*\otimes \mathbb{T}^*)+\mathrm{sym}(C_{\mathbb{T}}(\mathbb{T},.,.))\bigl]
\end{split}
\]
is a conformal invariant of weight one on $M_\bundle{E}$, in this case $\mathrm{sym}(b)$ is defined by $\mathrm{sym}(b)(x,y):=b(x,y)+b(y,x)$. Obviously $B_{\mathbb{T}}$ is the well known Bach tensor if $n=4$. Since $B_{\mathbb{T}}$ and $C_{\mathbb{T}}$ vanish for Einstein spaces, $B_{\mathbb{T}}=0$ as well as $C_{\mathbb{T}}=0$ are necessary conditions for a conformal Einstein space. Nevertheless, the class of Riemannian manifolds with $B_{\mathbb{T}}=0$ and $C_{\mathbb{T}}=0$ is larger than the class of conformal Einstein spaces (cf.~\cite{KNT,NP,Schmidt}).
\section{Preliminaries}

Let \( (M^{n},g) \) be a Riemannian manifold, then $\nabla $ denotes the Levi--Civita connection as well as \( R \) the Riemannian curvature tensor of $g$:
\[
R(X,Y,Z,T)=g(R_{X,Y}Z,T)=g(\nabla _{X}\nabla _{Y}Z-\nabla _{Y}\nabla _{X}Z-\nabla _{[X,Y]}Z,T)\, .\]
The Ricci tensor \( \ric \) is given by \( \ric(X,Y)=\mathrm{trace}\{V\mapsto R_{V,X}Y\} \) and the scalar curvature by \( S=\mathrm{trace}(\ric )\, . \)
Using the Kulkarni--Nomizu product:
\[ 
\begin{split}
(g\odot h)(X,Y,Z,T)  :=  g(X,T)h(Y,Z)+g(Y,Z)h(X,T)\qquad \\
   -g(X,Z)h(Y,T)-g(Y,T)h(X,Z)
\end{split} \]
we obtain the Weyl tensor \( W \) and the Schouten tensor \( \schout \)
:
\begin{equation}
\label{schout}
W:=R-g\odot \schout\quad ,\qquad \schout:=\frac{1}{n-2}\left( \ric-\frac{1}{2(n-1)}Sg\right) \, .
\end{equation}

Two Riemannian manifolds \( (M,g) \) and \( (N,h) \) are said to be \emph{conformally equivalent} if there are a diffeomorphism \( f:M\rightarrow N \) and a smooth function \( \psi :M\rightarrow (0,\infty ) \) satisfying \( f^{*}h=\psi ^{-2}g \). Since this diffeomorphism \( f \) is an isometry from \( (M,\overline{g}:=\psi ^{-2}g) \) to \( (N,h) \), we consider conformal transformations of the type: \( (M,g)\rightarrow (M,\overline{g}:=\psi ^{-2}g) \). The corresponding symbols for \( (M,\overline{g}) \) will be denoted by \( \overline{\nabla },\overline{R},\overline{W},... \) If \( (M,g)\rightarrow (M,\overline{g}:=\psi ^{-2}g) \) is a conformal transformation with \( \psi =e^{\phi } \) and \( \phi :M\rightarrow \mathbbm {R} \) smooth, then the Levi--Civita connections and the Weyl tensors are related by:
\begin{eqnarray}
\label{konfnabla}&
\overline{\nabla }_{X}Y=\nabla _{X}Y-\mathrm{d}\phi (X)Y-\mathrm{d}\phi (Y)X+\left\langle X,Y\right\rangle \nabla ^g\phi \, ,&\\
\label{konfweyl}
&\overline{W}=\psi ^{-2}W\, .&
\end{eqnarray}
In this case $\nabla ^g\phi $ is the vector field dual to $\mathrm{d}\phi $, i.e.~$(\nabla ^g\phi )^*=\mathrm{d}\phi $ while $*:TM\to T^*M$ is the isomorphism given by $Y^*(X)=g(X,Y)$.

Suppose $(\bundle{V},h)$ is a Riemannian vector bundle over $M$. If $A:\bundle{V}\to \bundle {V}$ is a symmetric bundle endomorphism w.r.t.~the inner product $h$, then at each point $p\in M$ there is an endomorphism $A^{\pen}_p\in \mathrm{End}(\bundle{V}_p)$ satisfying (Moore--Penrose inverse, cf.~\cite{Golan})
\begin{equation}
\begin{split}
\label{ident0}
&A_p\circ A_p^\pen \circ A_p=A_p\ , \qquad  A_p^\pen \circ A_p\circ A_p^\pen =A_p^\pen ,\\
&A_p \circ A_p^\pen  \ \ \text{and}\ \   \ A_p^\pen \circ A_p \ \ \text{are symmetric w.r.t.}\ h.
\end{split}
\end{equation}
Let $M_A$ consist of all points $p\in M$ at which the number of distinct eigenvalues of $A$ is locally constant. In particular, $M_A$ is open and dense in $M$. If $A\in \Gamma (\mathrm{End}(\bundle{V}))$ is of order $C^k$, the map
\[
A^\pen :M_A\to \mathrm{End}(\bundle{V}_{|M_A}), p\mapsto A^\pen _p\]
is of order $C^k$, i.e.~$A^\pen \in \Gamma (\mathrm{End}(\bundle{V}_{|M_A}))$. $A^\pen $ is called the \emph{Moore--Penrose inverse} of $A$.  In every point $p\in M$, $\bundle{V}$ admits an orthogonal decomposition into eigenspaces of $A$: $\bundle{V}_p=\bundle{V}^1_p\oplus \cdots \oplus \bundle{V}^k_p$. Since in the connected components of $M_A$, the number of distinct eigenvalues is constant, $p\mapsto \bundle{V}^j_p$ is differentiable, and thus $\bundle{V}^j\to U$ is a subbundle of $\bundle{V}\to U$ as long as $U\subset M$ is contained in a connected component of $M_A$. The existence and the uniqueness of the Moore--Penrose inverse follows from linear algebra, while the differentiability of $A^\pen $ goes as follows. Let $U$ be a connected component of $M_A$, then $\mathrm{Im}(A)\to U$ as well as $\mathrm{Ker}(A)\to U$ are subbundles of $\bundle{V}\to U$. In particular $A_{|U}$ supplies an orthogonal decomposition $\bundle{V}=\mathrm{Im}(A)\oplus \mathrm{Ker}(A)$. Now $A$ restricted to $\mathrm{Im}(A)$ is invertible, i.e.~set $A^\pen :=A^{-1}$ on $\mathrm{Im}(A)$ as well as $A^\pen =0$ on $\mathrm{Ker}(A)$. Since $\mathrm{Im}(A)$ and $\mathrm{Ker}(A)$ are subbundles of $\bundle{V}_{|U}$ as well as the assignment $A\mapsto A^{-1}$ is smooth, $A^\pen $ is on $M_A$ of order $C^k$ as long as $A$ is of order $C^k$ and moreover, $A^\pen $ satisfies the identities in (\ref{ident0}). Furthermore, if $\det A$ nowhere vanishes on $M$, $A^\pen $ equals $A^{-1}$ and is defined on all of $M$.

\begin{defn}
$\mathfrak{T}^r_s(M)$ may denote the set of $(r,s)$ tensor fields on a Riemannian manifold $(M,g)$. Let \( A\in \mathfrak {T}^{0}_{4}(M) \) be an algebraic curvature tensor. \( A \) becomes an endomorphism of \( \mathfrak {T}^{0}_{2}(M) \) in the following way:
\[
\curv{A}:\mathfrak {T}^{0}_{2}(M)\rightarrow \mathfrak {T}^{0}_{2}(M)\, ,\, X^*\otimes Y^*\mapsto \curv{A}(X^*\otimes Y^*) \]
where $ \curv{A}(X^*\otimes Y^*)(Z,T):=A(Y,Z,T,X)$. If $b\in \mathfrak{T}_2^0(M)$ is (skew) symmetric, $\curv{A}(b)$ is (skew) symmetric. In particular $\curv{A}$ is an endomorphism on two forms:
\[
\curv{A}:\Lambda ^2(T^*M)\to \Lambda ^2(T^*M).
\]
Since the first Bianchi identity implies
\[
g( \curv{A}(X^*\wedge Y^*),Z^*\wedge T^*) =\curv{A}(X^*\wedge Y^*)(Z,T)=A(X,Y,Z,T),
\]
$\curv{A}$ is symmetric on $\Lambda ^2(T^*M)$ w.r.t.~the extension of $g$. In this case $X^*\wedge Y^*$ is the two form given by $X^*\otimes Y^*-Y^*\otimes X^*$.
\end{defn}
Let $W$ be the Weyl tensor on $(M^n,g)$ and define
\[
\bundle{E}:=\bigcup _{p \in M}\bundle{E}_p, \ \ \text{while}\ \ \bundle{E}_p:=\{ v\in T_pM| W(v,.,.,.)=0 \} .
\]
$\bundle{E}_p$ are vector spaces, but in general $\bundle{E}$ is not a vector bundle of $M$. If $\bundle{E}^\perp $ is the orthogonal complement of $\bundle{E}$ in $TM$, the decomposition $\Lambda ^2(TM)=\mathrm{Ker}(\curv{W})\oplus \mathrm{Im}(\curv{W})$ as well as
\[
\frak{so}(T_pM)=\frak{so}(\bundle{E}_p^\perp )\oplus \frak{so}(\bundle{E}_p)\oplus \bundle{E}_p^\perp \otimes \bundle{E}_p
\] 
imply:
\begin{equation}
\label{gl10}
\dim \mathrm{Ker}(\curv{W}_p)\geq \dim \bundle{E}_p\left( n-\frac{\dim \bundle{E}_p+1}{2}\right) .
\end{equation}
Contrary a lower bound of $\dim \bundle{E}_p$ in terms of $\dim \mathrm{Ker}(\curv{W}_p)$ is not possible in general, since if $M$ is even--dimensional and there is a non--degenerate two form $\eta $ in the image of $\curv{W}_p$, $v\llcorner \eta \neq 0$ for all $v\in T_pM-\{ 0\} $ implies $\dim \bundle{E}_p=0$. 
 
Denote by $\frak{w}: T^*M\to T^*M$ the negative Ricci contraction of $ \curv{W}^2=\curv{W}\circ \curv{W}$, i.e.~$\frak{w}(\Theta )(X)$ is the trace of the mapping
\[
T^*M\to T^*M,\ \eta \mapsto X\llcorner \curv{W}^2(\theta \wedge \eta ) \ .
\]
If $f_1,\ldots ,f_n$ is a base of $T_pM$ and $\eta _1,\ldots ,\eta _n$ the corresponding cobase of $T^*_pM$ ($\eta _j(f_i)=\delta _{ij}$), we have
\[
\frak{w}(\theta )=\sum _{i=1}^nf_i\llcorner \curv{W}^2(\eta _i\wedge \theta ).
\]
Using an orthonormal base w.r.t.~$g$, we conclude that $\frak{w}$ is symmetric w.r.t.~$g$ and in particular, $ \curv{W}^2\geq 0$ supplies $\frak{w}\geq 0$. Since $\frak{w}$ is symmetric, there is an open and dense subset $M_\frak{w}$ of $M$ on which the Moore--Penrose inverse of $\frak{w}$ exists. Moreover, $*:TM\to T^*M, v\mapsto g(v,.)$ gives an isomorphism
\begin{equation}
\label{gl12}
\bundle{E}_p\to \mathrm{Ker}(\frak{w}_p).
\end{equation}
The fact $*(\bundle{E}_p )\subseteq \mathrm{Ker}(\frak{w}_p) $ is obvious from the definitions. In order to see equality, let $v\in T_pM$ be a vector with $\frak{w}(v^*)=0$. We have to show $\curv{W}(v^*\wedge \theta )=0$ or equivalently $\curv{W}^2(v^*\wedge \theta )=0$ for all $\theta \in T^*_pM$. Suppose $\eta _1,\ldots ,\eta _m\in \Lambda ^2(T^*_pM)$ is an orthonormal base of eigenvectors to the non--negative eigenvalues $\lambda _1,\ldots ,\lambda _m$ of $\curv{W}^2$, then $\curv{W}^2_p$ considered as $(0,4)$ tensor is given by $\sum _j\lambda _j\eta _j\otimes \eta _j$. Thus, considering the two forms $\eta _j$ as skew symmetric maps $T^*_pM\to T^*_pM$ leads to $\frak{w}_p=-\sum _{j=1}^m\lambda _j(\eta _j)^2$. Since each of these summands is non--negative definite, $\frak{w}(v^*)=0$ implies $(\eta _j)^2(v^*)=0$ for all $j$ with $\lambda _j\neq 0$. But this gives the claim $\curv{W}^2(v^*\wedge \theta )=0$. Therefore, equation (\ref{gl12}) and the above arguments supply that $\bundle{E}\to M_\frak{w}$ is smooth, in particular, if $U$ is a connected component of $M_\frak{w}$, $\bundle{E}_{|U}$ is a subbundle of $TM_{|U}$. In the introduction we defined $M_\bundle{E}$ to be set of all points $p\in M$ at which the rank of $\bundle{E}$ is locally constant, in particular $M_\frak{w}\subseteq M_\bundle{E}$. Considering the map $\frak{w}^\pen \frak{w}$ as endomorphism on $TM_\frak{w}$ (using the $*$ isomorphism), $1-\frak{w}^\pen \frak{w}$ is the projection $TM_\frak{w}\to \bundle{E}$, i.e.~$1-\frak{w}^\pen \frak{w}$ can be differentiable extended to $M_\bundle{E}$. Using (\ref{gl12}) and the fact $\frak{w}^\pen =(\frak{w}_{|\mathrm{Im}(\frak{w})})^{-1}$, $\frak{w}^\pen $ is differentiable on $M_\bundle{E}$. Thus, we can assume $M_\frak{w}=M_\bundle{E}$.

Since the Weyl tensor is conformally invariant, it should be mentioned that $\bundle{E}$, $\mathrm{Ker}(\curv{W})$ as well as $M_\bundle{E}$ are invariant under conformal transformation.
\begin{rem}
\label{rema1}
Let $(M,g)$ be a four dimensional Riemannian manifold. Suppose $U$ is a connected component of $M_\bundle{E}$, then either $(U,g)$ is conformally flat (i.e.~$W_p=0$ for all $p\in U$) or the rank of the bundle $\bundle{E}\to U$ is zero. In order to see this suppose $0\neq v\in \bundle{E}$ and $|v|=1$. Let $v,e_1,e_2,e_3$ be an orthonormal base of $T_pM$, then $v^*\wedge e_j^*$ and $\hodge (v^*\wedge e_j^*)$, $j=1,2,3$, give a base of $\Lambda ^2(T_p^*M)$. Since $\curv{W}$ commutes with the Hodge star operator $\hodge $, we obtain $\curv{W}(\eta )=0$ for any two form $\eta $. But this implies $W=0$. Moreover, the Ricci contraction of $\curv{W}^2$ satisfies (cf.~\cite[16.75]{Bes} resp.~\cite{De1})
\[
\frak{w}=\frac{1}{2}\mathrm{tr}(\curv{W}^2)\mathrm{Id}.
\]
\end{rem}

\section{Spaces of harmonic Weyl tensor}
Let $(M,g)$ be a Riemannian manifold and $A$ be a $(0,s+1)$ tensor, then define for $r\leq s+1$ the $(0,s)$ tensor field
\[
\delta _rA(X_1,\ldots ,X_s):=\sum _{j=1}^n(\nabla _{e_j}A)(X_1,\ldots ,X_{r-1},e_j,X_r,\ldots X_s)
\]
($e_1,\ldots ,e_n$ is an orthonormal frame). $\delta _rA$ is called the divergence of $A$ and if $A$ is symmetric, $\delta _rA$ does not depend on $r$. If $A$ is a curvature operator, introduce as abbreviation $\delta A:=\delta _4A$.
Moreover, we define the exterior derivative of a symmetric $(0,2)$ tensor $b$ by
\[
\mathrm{d} ^\snabla  b(X,Y,Z):=(\nabla _Xb)(Y,Z)-(\nabla _Zb)(X,Z).
\]
The differential Bianchi identity supplies (cf.~\cite[Ch.~16.3]{Bes})
\begin{equation}
\label{deltaW}
\delta R=\mathrm{d} ^\snabla \mathrm{Ric},\qquad \delta W=(n-3)\mathrm{d} ^\snabla  \schout
\end{equation}
where $\schout$ is the Schouten tensor given in (\ref{schout}). A Riemannian manifold \( (M^{n},g) \) of dimension \( n\geq 4 \) is called \emph{C-space} or \emph{space with harmonic Weyl tensor} if $\delta W=0$ (cf. \cite[(16.D)]{Bes}). Let $\overline{g}=\psi ^{-2}g$ with $\psi =e^\phi $ be a conformal transformation, we obtain the following relation:
\begin{equation}
\label{gl52}
\overline{\delta } \overline{W}=\delta W-(n-3)W(.,.,.,\nabla ^g\phi ).
\end{equation}
In order to see this, use equations (\ref{konfweyl}) and (\ref{konfnabla}) or cf.~\cite[16.25]{Bes}. In particular a Riemannian manifold $(M^n,g)$, $n\geq 4$, is conformally related to a space with harmonic Weyl tensor if and only if there is a function $\phi $ with
\begin{equation}
\label{gl20}
\delta W_Z=(n-3)\curv{W}(Z^*\wedge \mathrm{d}\phi ),
\end{equation}
in this case $\delta W_Z$ is the two form given by
\[
\delta W_Z(X,Y):=\delta W(X,Y,Z).
\]
A necessary condition for a solution of (\ref{gl20}) is $\delta W_Z\in \mathrm{Im}(\curv{W})$ or equivalent $\curv{W}\curv{W}^\pen (\delta W_Z)=\delta W_Z$, but we compute $\mathrm{d}\phi $ explicitly, i.e.~this condition will be superfluous. Let $\frak{w}$ be the negative Ricci contraction of $\curv{W}^2$ and $\frak{w}^\pen $ be the Moore--Penrose inverse of $\frak{w}$ which is well defined and smooth on $M_\bundle{E}$. Apply $\curv{W}$ to equation (\ref{gl20}) to obtain
\[
\curv{W}(\delta W_Z)=(n-3)\curv{W}^2 (Z^*\wedge \mathrm{d}\phi ).
\]
Moreover, suppose $e_1,\ldots ,e_n$ is a $g$--orthonormal base of $T_pM$, then $e_1^*,\ldots ,e_n^*$ is the cobase in $T_p^*M$ and we get
\[
\sum _{i=1}^ne_i\llcorner \curv{W}(\delta W_{e_i})=(n-3)\sum _{i=1}^ne_i\llcorner \curv{W}^2(e_i^*\wedge \mathrm{d}\phi )=(n-3)\frak{w}(\mathrm{d}\phi ) .
\]
Thus, on the components of $M_\bundle{E}$ where $\mathrm{rank}(\bundle{E})=0$, we obtain
\[
\mathrm{d}\phi = \frac{1}{n-3}\frak{w}^\pen  \left( \sum _{i=1}^n e_i\llcorner \curv{W} (\delta W_{e_i}) \right) .
\]
\begin{defn}
Let $(M^n,g)$ be a smooth Riemannian manifold of dimension $n\geq 4$ and $e_1,\ldots ,e_n$ be a $g$--orthonormal base, then the vector field $\mathbb{T}$ given by
\begin{equation}
\label{definition_of_T}
g(\mathbb{T},.):=\frac{1}{n-3}\frak{w}^\pen \left( \sum _{i=1}^n e_i\llcorner \curv{W}(\delta W_{e_i}) \right)
\end{equation}
as well as the $(0,3)$ tensor
\[
C_\mathbb{T}:=\mathrm{d} ^\snabla \, \schout-W(.,.,.,\mathbb{T})
\]
are smooth on the open and dense subset $M_\bundle{E}\subseteq M$.
\end{defn}
\begin{rem}
Let $(M,g)$ be a Riemannian manifold of dimension $4$, then $M_\bundle{E}=M_0\cup M_1$, where $M_0$ is the interior of $\{ p\in M| \ W(p)=0\}$ and as well as $M_1$ is the open set $\{ p\in M| \ W(p)\neq 0\} $. In particular $\mathbb{T}$ vanishes on $M_0$ and equals 
\[
\frac{2}{\mathrm{tr}(\curv{W}^2)}\sum _{i,j=1}^4\curv{W}(\delta W_{e_i})(e_i,e_j)e_j
\]
on $M_1$.
\end{rem}
In the definition of $C_\mathbb{T}$ we used $\mathrm{d} ^\snabla \schout $ instead of $\delta W$ in order to get the Cotton tensor in dimension $n=3$, but for the computations we consider
\[
C_\mathbb{T}=\frac{1}{n-3}\delta W-W(.,.,.,\mathbb{T}).
\]
\begin{prop}
Suppose $\overline{g}=\psi ^{-2}g$ is a conformal transformation with $\psi =e^\phi$, then $\overline{C} _{\overline{\mathbb{T}}}= C_\mathbb{T}$ holds on $M_\bundle{E}$.
\end{prop}
\begin{proof}
$\overline{W}=\psi ^{-2}W$ supplies $\overline{\curv{W}}=\psi ^2\curv{W}$, i.e.~we obtain $\overline{\frak{w}}=\psi ^4\frak{w}$ as well as $\overline{\frak{w}}^\pen =\psi ^{-4}\frak{w}^\pen $. If $e_1,\ldots ,e_n$ is an orthonormal base w.r.t.~$g$, $\overline{e_j}:=\psi e_j$ gives an orthonormal base w.r.t.~$\overline{g}$, i.e.
\begin{eqnarray*}
\overline{\mathbb{T}}^{\overline{*}}&=&\frac{1}{n-3}\overline{\frak{w}}^\pen \left(\sum _{i=1}^n\overline {e_i}\llcorner \overline{\curv{W}}^t (\overline{\delta }\overline{W}_{\overline{e_i}}) \right) \\
&=& \frac{1}{n-3}\frak{w}^\pen \left( \sum _{i=1}^n e_i \llcorner \curv{W}^t \left( \delta W_{e_i}-(n-3)\curv{W}(e_i^*\wedge \mathrm{d}\phi )\right) \right) \\
&=& \mathbb{T}^*-\frak{w}^\pen \frak{w}(\mathrm{d}\phi ).
\end{eqnarray*}
Thus, the vector fields are related by
\begin{equation}
\label{gl30}
\overline{\mathbb{T}}=\psi ^2\mathbb{T}-\psi ^2(\nabla ^g\phi )_{\bundle{E}^\perp},
\end{equation}
where $(\nabla ^g\phi )_{\bundle{E}^\perp}$ is the projection of $\nabla \phi $ to the orthogonal complement of $\bundle{E}$ in $TM$. Since $W(.,.,.,(\nabla \phi )_{\bundle{E}})=0$, we conclude the conformal invariance of the tensor $C_\mathbb{T}$ if $n>3$ [use (\ref{deltaW}) and (\ref{gl52})]:
\begin{eqnarray*}
\overline{C}_{\overline{\mathbb{T}}}&=&\frac{1}{n-3}\overline{\delta }\overline{W}-\overline{W}(.,.,.,\overline{\mathbb{T}})\\
&=&\frac{1}{n-3}\delta W-W(.,.,.,\nabla ^g\phi )-W(.,.,.,\mathbb{T}-(\nabla ^g\phi )_{\bundle{E}^\perp})\ = \ C_\mathbb{T}.
\end{eqnarray*}
If $n=3$, $C_\mathbb{T}$ is the Cotton tensor whose conformal invariance is already known in this dimension.
\end{proof}

Since $C_\mathbb{T}$ is trivial on a space with harmonic Weyl tensor, the vanishing of $C_\mathbb{T}$ is a necessary condition for a Riemannian manifold to be locally conformally related to a C--space. However, $C_\mathbb{T}=0$ is a sufficient condition if and only if there is a section $V$ in $\bundle{E}$ in such a way that $\mathbb{T}+V$ is a (differentiable) gradient field on $M$.
\begin{prop}
Let $(M^n,g)$ be a simply connected Riemannian manifold of dimension $n\geq 4 $ such that $\mathrm{rank}(\bundle{E})=0$ on $M_\bundle{E}$. Then $(M,g)$ is (globally) conformally related to a space with harmonic Weyl tensor if and only if $\mathrm{d}\mathbb{T}^*=0$ and $C_\mathbb{T}=0$ on $M_\bundle{E}$ as well as $\mathbb{T}$ can be (differentiable) extended to $M$. Moreover, if $\det \curv{W}$ does not vanish on $M_\bundle{E}$, the condition $\mathrm{d}\mathbb{T}^*=0$ follows from $C_\mathbb{T}=0$.
\end{prop}
\begin{proof}
The claims are obvious from the uniqueness of the conformal factor (up to scaling) as well as the definition of $\mathbb{T}$ and $C_\mathbb{T}$. Since $M$ is simply connected and $\mathbb{T}^*$ is exact, there is a function $\phi :M\to \mathbbm{R}$ with $\nabla ^g\phi =\mathbb{T}$. Set $\psi :=e^\phi $, then $\psi ^{-2}g$ is a space with harmonic Weyl tensor. In order to see the last claim consider the divergence of $C_\mathbb{T}$ w.r.t.~the third argument. A straightforward calculation shows (cf.~\cite{List3})
\[
\delta _3(C_\mathbb{T})=C_\mathbb{T}(.,.,\mathbb{T})-\curv{W}(\mathrm{d}\mathbb{T}^*) ,
\]
i.e.~we conclude $\mathrm{d}\mathbb{T}^*=0$ from the injectivity of $\curv{W}$ and $C_\mathbb{T}=0$.
\end{proof}

\section{Conformal Einstein Spaces}

A Riemannian manifold \( (M,g) \) of dimension \( n\geq 3 \) is called
\emph{Einstein space} if the traceless Ricci tensor $\mathrm{Ric}^{\circ }:=\mathrm{Ric}-\frac{S}{n}g$ vanishes. In this case, \( S \) is constant. $(M,g)$ is said to be a \emph{conformal Einstein space} if $g$ is locally conformally related to an Einstein space. Conformal Einstein spaces have already been considered in the 1920's by Brinkmann (cf.~\cite{Brink1}). 

Let \( (M,g)\rightarrow (M,\overline{g}:=\psi ^{-2}g) \) with \( \psi =e^{\phi } \) be a conformal transformation. The Ricci tensor has the following transformation behavior (cf. \cite[Lemma A.1]{Kh} or \cite{Bes}):
\begin{equation}
\label{konfRic}
\overline{\mathrm{Ric}}=\mathrm{Ric}+(n-2)[\nabla ^{2}\phi +\mathrm{d}\phi \otimes \mathrm{d}\phi ]+\left[ \triangle \phi -(n-2)\left\langle \nabla ^g\phi ,\nabla ^g\phi \right\rangle \right] g\, ,
\end{equation}
where \( \nabla ^{2}\phi  \) is the Hessian of \( \phi  \) (i.e. \( \nabla ^{2}\phi (X,Y)=\left\langle \nabla _{X}\nabla ^g\phi ,Y\right\rangle  \))
and $\triangle \phi $ is the trace of $\nabla ^2\phi $. If \( (M,\overline{g}) \)
is an Einstein space, we conclude taking the trace of the last equation for
the Ricci tensor of \( (M,g) \) :
\begin{equation}
\label{ricci}
0=\mathrm{Ric}^\circ +(n-2)[\nabla ^{2}\phi +\mathrm{d}\phi \otimes \mathrm{d}\phi ]-\frac{n-2}{n}\left[ \triangle \phi +\left| \nabla ^g\phi \right| ^{2}\right] g\, .
\end{equation}

\begin{defn}
Let \( V \) be a vector field, then the traceless \( (0,2) \) tensor field \( E_{V} \) given by
\[
E_V:=\mathrm{Ric}^\circ +(n-2)\Bigl(\nabla V^*+V^*\otimes V^*-\frac{1}{n}\left[ \mathrm{div}(V)+g(V,V)\right]g\Bigl)
\]
is called \emph{conformal Ricci tensor w.r.t.~$V$}.
\end{defn}
The tensor field $E_V$ is symmetric if and only if $\nabla V^*$ is symmetric, i.e.~if and only if $V^*$ is closed. Thus a Riemannian manifold $(M,g)$ is locally conformally related to an Einstein space if and only if there is a vector field $V$ with $E_V=0$. If additionally $M$ is simply connected and $E_V$ vanishes, there is a function $\psi =e^\phi $ which gives the Einstein space $(M,\psi ^{-2}g)$. In this case $V$ equals $\nabla ^g\phi $. The differential Bianchi identity tells that an Einstein space has a harmonic Weyl tensor [cf.~equation (\ref{deltaW})]. Thus, the first candidate of a vector field such that $E_V=0$ holds, is the vector field $\mathbb{T}$ given in (\ref{definition_of_T}).

\begin{lem}
Let $\overline{g}=\psi^{-2 }g$, $\psi =e^\phi $, be a conformal transformation and $\overline{\mathbb{T}}$ as well as $\mathbb{T}$ be the corresponding vector fields for $\overline {g}$ and $g$ defined in (\ref{definition_of_T}). If $\nabla ^g\phi $ is contained in $\bundle{E}^\perp \subset TM$ in every point of $M_\bundle{E}$, the conformal Ricci tensors satisfy:
\[
\overline{E}_{\overline{\mathbb{T}}}=E_\mathbb{T}.
\]
\end{lem}
\begin{proof}
Introduce the $(0,2)$ tensor
\[
F_V:=\nabla V^*+V^*\otimes V^*-\frac{1}{n}\left[ \mathrm{div}(V)+g(V,V)\right]g . 
\]
One easily verifies for two vector fields $V$ and $Z$:
\[
F_{V+Z}=F_V+F_Z+V^*\otimes Z^*+Z^*\otimes V^*-\frac{2}{n}g(V,Z)g.
\]
Moreover, using (\ref{konfnabla}) a straightforward calculation yields
\[
\overline{F}_{\psi ^2V}= F_V+\mathrm{d}\phi \otimes V^*+V^*\otimes \mathrm{d}\phi -\frac{2}{n}g(\nabla ^g\phi ,V)g.
\]
Since $\psi ^{-2}$ is the  conformal factor, equation (\ref{gl30}) yields $\overline{\mathbb{T}}=\psi ^2(\mathbb{T}-Y)$, where $\nabla ^g\phi =X+Y $ is the orthogonal decomposition in $(\bundle{E}\oplus \bundle {E}^\perp )$, i.e.~$X:M_\bundle{E}\to \bundle{E}$ and $Y:M_\bundle{E}\to \bundle{E}^\perp $ are sections and $\bundle{E}^\perp $ is the $g$--orthogonal complement of $\bundle{E}$ in $TM$. Thus, we obtain
\begin{eqnarray*}
\overline{F}_{\overline{\mathbb{T}}}&=&F_{\mathbb{T}-Y }+\mathrm{d}\phi \otimes (\mathbb{T}^*-Y^*)+(\mathbb{T}^*-Y^*)\otimes \mathrm{d}\phi -\frac{2}{n}g(\nabla ^g\phi ,\mathbb{T}-Y )g\\
&=&F_\mathbb{T}+F_{-Y }+X^*\otimes \mathbb{T}^*+\mathbb{T}^*\otimes X^*-2Y^*\otimes Y^*+\frac{2}{n}g(Y,Y)g\\
&&-X^*\otimes Y^*-Y^*\otimes X-\frac{2}{n}g(X,\mathbb{T}-Y)g\\
&=&F_\mathbb{T}-F_{Y }+X^*\otimes (\mathbb{T}^*-Y^*)+(\mathbb{T}^*-Y^*)\otimes X^*-\frac{2}{n}g(X,\mathbb{T}-Y)g.
\end{eqnarray*}
Use equation (\ref{konfRic}) to show that the traceless Ricci tensors are related by $(n-2)F_{\nabla \phi}$, i.e.~the definition of $E_V$ gives:
\begin{equation}
\begin{split}
\overline{E}_{\overline{\mathbb{T}}}&=\overline{\mathrm{Ric}}^\circ +(n-2)\overline{F}_{\overline{\mathbb{T}}}\\
&=\mathrm{Ric}^\circ +(n-2)\Bigl[ F_{X+Y }+F_\mathbb{T}-F_Y+\\
&\quad +X^*\otimes (\mathbb{T}^*-Y^*)+(\mathbb{T}^*-Y^*)\otimes X^* -\frac{2}{n}g(X,\mathbb{T}-Y)\Bigl] \\
\label{konE}
&=E_\mathbb{T}+(n-2)\Bigl[ F_{X}+X^*\otimes \mathbb{T}^*+\mathbb{T}^*\otimes X^*-\frac{2}{n}g(X,\mathbb{T})g\Bigl] .
\end{split}
\end{equation}
But we assumed $\nabla ^g\phi \in \bundle{E}^\perp $, i.e.~$X=0$ supplies the claim.
\end{proof}
\begin{rem}
Thus, we have proved that $E_\mathbb{T}$ is a conformal invariant on $M_\bundle{E}$ for all Riemannian manifolds $(M,g)$ with $\mathrm{rank}(\bundle{E})=0$ on $M_\bundle{E}$. In particular, if $(M,g)$ is a Riemannian four manifold, $E_\mathbb{T}$ is conformally invariant on the open subset
\[
\{ p\in M\ |\ W(p)\neq 0\}\subseteq M.\]
\end{rem}
\begin{thm}
\label{main_thm}
Let $(M,g)$ be a simply connected Riemannian manifold of dimension $n\geq 4$ such that $\mathrm{rank}(\bundle{E})=0$ on $M_\bundle{E}$. Then $(M,g)$ is conformally related to an Einstein space if and only if $E_\mathbb{T}$ vanishes on $M_\bundle{E}$ and $\mathbb{T}$ is extendible to a $C^2$ vector field on $M$. 
\end{thm}
\begin{proof}
Since $E_\mathbb{T}$ is a conformal invariant on manifolds with $\mathrm{rank}(\bundle{E})=0$ and $\mathbb{T}$ vanishes on Einstein spaces, $E_\mathbb{T}=0$ is a necessary condition. Conversely, $E_\mathbb{T}=0$ implies that $\mathbb{T}^*$ is closed, i.e.~there is a $C^3$ function $\phi :M\to \mathbbm{R}$ with $\mathbb{T}=\nabla \phi $. Then the above computations show that $e^{-2\phi }g$ is an Einstein metric on $M$ with $C^1$ Ricci tensor, and therefore analytic Ricci tensor in suitable coordinates (cf.~\cite{DeT}).
\end{proof}
\begin{rem}
If $(M,g)$ is a connected Einstein space, the Weyl tensor is real analytic in suitable coordinates, i.e.~$W$ vanishes on an open subset of $M$ if and only if $g$ is of constant sectional curvature. Therefore, a Riemannian $4$--manifold which is conformally related to an Einstein space is conformally flat or $\mathrm{rank}(\bundle{E})=0$ on $M_\bundle{E}$. Apply the above theorem which was already mentioned in case $\dim M=4$ in \cite{List3}, then a Riemannian manifold $(M^4,g)$ is (locally) conformally related to an Einstein space if and only if $g$ is conformally flat or $\mathbb{T}$ is extendible to a vector field on $M$ and $E_\mathbb{T}$ vanishes.
\end{rem}

\section{Generalized Bach tensor}
In four dimensions there is another conformal invariant called \emph{Bach tensor} (cf.~\cite{Bach}). The Bach tensor is defined by
\[
B:=\delta _1\delta _4W+\frac{1}{2}\curv{W}(\mathrm{Ric}).
\]
If $(M,g)$ is locally conformally related to an Einstein space, $B$ has to vanish, but there are Riemannian manifolds with $B=0$ which are not conformally related to Einstein spaces (cf.~\cite{Schmidt}). If $M$ is a compact $4$--manifold, the Bach tensor is the gradient of the functional
\[
g\mapsto \int\limits _M|W|^2\mathrm{vol}_g,
\]
in particular, $B$ vanishes for critical metrics of this functional.

Let $\overline{g}=\psi ^{-2}g$ with $\psi =e^\phi $ be a conformal transformation on $M^n$. Since $C_V$ satisfies the first Bianchi identity and $\sum C_V(e_i,.,e_i)=0$, we obtain
\[
\begin{split}
\overline{\delta }_1\overline{\delta }_4\overline{W}(x,y)=&(n-3)\overline{\delta }_1C_{\nabla \phi }(x,y)\\
=& (n-3)\psi ^2\bigl[ \delta _1C_{\nabla \phi }(x,y)-(n-4)C_{\nabla \phi }(\nabla \phi ,x,y)+\\&\qquad \qquad +C_{\nabla \phi }(\nabla \phi ,y,x)\bigl] .
\end{split} \]
Moreover, $\delta _1W(X,Y,Z)=\delta _4W(Z,Y,X)$ implies
\[
\begin{split}
(n-3)\delta _1C_{\nabla \phi }(x,y)=&\delta _1\delta _4W(x,y)-(n-3)\delta _1W(x,y,\nabla \phi )\\
&\ -(n-3)\curv{W}(\nabla ^2\phi )(x,y)\\
=&\delta _1\delta _4W(x,y)-(n-3)^2C_{\nabla \phi }(\nabla \phi ,y,x)\\
&-(n-3)^2\curv{W}(\mathrm{d}\phi \otimes \mathrm{d}\phi )(x,y)-(n-3)\curv{W}(\nabla ^2\phi )(x,y),
\end{split}\]
i.e.~if $\mathrm{sym}(b)(x,y):=b(x,y)+b(y,x)$ denotes the symmetrization of a $(0,2)$ tensor, this leads to
\begin{equation}
\begin{split}
\psi ^{-2}\overline{\delta }_1\overline{\delta }_4\overline{W}=&\ \delta _1\delta _4W-(n-3)\curv{W}(\nabla ^2\phi +\mathrm{d}\phi \otimes \mathrm{d}\phi )\\
\label{gl59}
&-(n-4)(n-3)\left[ \mathrm{sym}(C_{\nabla \phi }(\nabla \phi ,.,.))+\curv{W}(\mathrm{d}\phi \otimes \mathrm{d}\phi )\right] .
\end{split}
\end{equation}
\begin{defn}
Let $(M,g)$ be a smooth Riemannian manifold of dimension $n\geq 4$ and $\mathbb{T}$ be the vector field given in (\ref{definition_of_T}). Then the trace--free, symmetric $(0,2)$ tensor
\[ \begin{split}
B_\mathbb{T}:=&\ \delta _1\delta _4W+\frac{n-3}{n-2}\curv{W}(\mathrm{Ric})-\frac{n-4}{n-2}\curv{W}(\mathrm{sym}(E_\mathbb{T}))\\
&\quad -(n-3)(n-4)\bigl[ \curv{W}(\mathbb{T}^*\otimes \mathbb{T}^*)+\mathrm{sym}(C_\mathbb{T}(\mathbb{T},.,.))\bigl]
\end{split} \]
is well defined and smooth on $M_\bundle{E}$ and will be called \emph{generalized Bach tensor}.
\end{defn}
\begin{prop}
The generalized Bach tensor is a conformal invariant on $M_\bundle{E}$ of weight $1$, i.e.~if $\overline{g}=\psi ^{-2}g$ is a conformal transformation, we have $\overline{B}_{\overline{\mathbb{T}}}=\psi ^2B_\mathbb{T}$.
\end{prop}
\begin{proof}
The Weyl tensor considered as endomorphism on $\frak{T}^0_2$ transforms like $\overline{\curv{W}}=\psi ^2\curv{W}$. Remember that $\overline{C}_{\overline{\mathbb{T}}}=C_\mathbb{T}$ as well as $\overline{\mathbb{T}}=\psi ^2(\mathbb{T}-\nabla ^g\phi _{\bundle{E}^\perp })$ hold, where $\bundle{E}\oplus \bundle{E}^\perp $ is the orthogonal decomposition w.r.t.~$g$. Since the Ricci contraction of $W$ vanishes, we have $\curv{W}(\lambda g)=0$ for any function $\lambda $, i.e.~taking $W(\nabla ^g\phi _{\bundle{E}},.,.,.)=0 $ as well as equations (\ref{gl59}), (\ref{konfRic}) and (\ref{konE}) into consideration leads to
\[ \begin{split}
\psi ^{-2}\overline{B}_{\overline{\mathbb{T}}}
=&\ B_\mathbb{T}-(n-4)\mathrm{sym}(\delta W(\nabla ^g\phi ,.,.))+(n-4)\mathrm{sym}(\delta W(\nabla ^g\phi _{\bundle{E}^\perp },.,.))\\
&\qquad \qquad \quad -(n-4)\curv{W}(\mathrm{sym}(\nabla (\mathrm{d} \phi _{\bundle{E} }) +\mathrm{d}\phi _{\bundle{E} }\otimes \mathrm{d} \phi _{\bundle{E}})).
\end{split}\]
$W(\nabla ^g\phi _\bundle{E},.,.,.)=0$ supplies for an orthonormal base of $T_pM$
\[
\delta W(\nabla ^g\phi _\bundle{E},.,.)=-\sum _iW(\nabla _{e_i}\nabla ^g\phi _\bundle{E},.,.,e_i),
\]
and since $\curv{W}$ commutes with $\mathrm{sym}$, we obtain
\begin{equation}
\label{gl63}
\mathrm{sym}(\delta W(\nabla ^g \phi _\bundle{E},.,.))=-\curv{W}(\mathrm{sym}(\nabla (\mathrm{d}\phi _\bundle{E}))).
\end{equation}
But this implies the claim: $\overline{B}_{\overline{\mathbb{T}}}:=\psi ^2B_\mathbb{T}$.
\end{proof}
\begin{lem}
Suppose $C_\mathbb{T}$ vanishes, then
\[
\widehat{B}_\mathbb{T}:=\delta _1\delta _4W+\frac{n-3}{n-2}\curv{W}(\mathrm{Ric})-(n-3)(n-4)W(\mathbb{T},.,.,\mathbb{T})
\]
and $\curv{W}(\mathrm{sym}(E_\mathbb{T}))$ are conformally invariant of weight $1$.
\end{lem}
\begin{proof}
In order to see the claim, we use equation (\ref{konE}) and the fact that in this equation $\nabla \phi _\bundle{E}=X\in \bundle{E}$:
\[
\overline{\curv{W}}\bigl( \mathrm{sym}(\overline{E}_{\overline{\mathbb{T}}})\bigl) =\psi ^2\curv{W}\bigl( \mathrm{sym}(E_\mathbb{T})+(n-2)\mathrm{sym}(\nabla (\mathrm{d} \phi _\bundle{E}))\bigl) .
\]
We apply equation (\ref{gl63}) to conclude
\[\begin{split}
\curv{W}(\mathrm{sym}(\nabla (\mathrm{d} \phi _\bundle{E})))=&\ -\mathrm{sym}(\delta W(\nabla ^g\phi _\bundle{E},.,.))\\
=&\ -(n-3)\mathrm{sym}(C_\mathbb{T}(\nabla ^g\phi _\bundle{E},.,.))=0.
\end{split}\]
The conformal invariance of $\widehat{B}_\mathbb{T}$ is then obvious from the definition of $B_\mathbb{T}$ and $C_\mathbb{T}=0$.
\end{proof}
If $(M,g)$ is an Einstein space, the conformal invariants $C_\mathbb{T}$ and $B_\mathbb{T}$ vanish. Thus, $C_\mathbb{T}$ and $B_\mathbb{T}$ are trivial for any Riemannian manifold which is locally conformally related to an Einstein space. Moreover, the last lemma tells that $\widehat{B}_\mathbb{T}$ and $\curv{W}(\mathrm{sym}(E_\mathbb{T}))$ have to vanish for a conformal Einstein space. However, $C_\mathbb{T}=0$, $\curv{W}(\mathrm{sym}(E_\mathbb{T}))=0$ and $\widehat{B}_\mathbb{T}=0$ does not seem to be sufficient for $(M,g)$ to be locally conformally related to an Einstein space (cf.~\cite{KNT,NP,Schmidt}).

\bibliographystyle{abbrv}
\bibliography{conv_inv.bbl}
\end{document}